\documentclass{amsart}

\usepackage{epsfig}
\usepackage{eufrak}
\usepackage{amsmath}
\newtheorem{Proposition}{Proposition}
  \newtheorem{Remark}{Remark}
  
  \newtheorem{Lemma}[Proposition]{Lemma}
  
  \newtheorem{Theorem}[Proposition]{Theorem}

\newtheorem{Definition}[Proposition]{Definition}

\def\la{\lambda}
\def\e{\epsilon}

\def\RR{\mathbb{R}}
\newcommand\ZZ{\mathbb{Z}}

\def\QED{\hfill\hfill\qed}

\def\z{\noindent}
\def\CC{\mathbb{C}}
\def\RR{\mathbb{R}}
\def\NN{\mathbb{N}}
\def\lap{\mathcal{L}}
\def\bor{\mathcal{B}}
\begin{document}
\begin{abstract}
  Borel summation techniques are developed to obtain exact invariants from
  formal adiabatic invariants (given as divergent series in a small parameter)
  for a class of differential equations, under assumptions of analyticity of
  the coefficients; the method relies on the study of associated partial
  differential equations in the complex plane.  The type and location of the
  singularities of these associated functions, important in determining
  exponentially small corrections to formal invariants are also briefly
  discussed.

\smallskip

\end{abstract}

\author{O. Costin, L. Dupaigne, and M. D. Kruskal}
\address{Department of Mathematics\\ Rutgers University
\\ Piscataway, NJ 08854-8019}

\title{Borel summation of adiabatic invariants}

\gdef\shorttitle{Borel summation of adiabatic invariants}
\maketitle

\section{
Introduction and main results}

For many ordinary differential equations depending on a small parameter (say
$\epsilon$) it is possible to construct adiabatic invariants: these are
generically divergent expansions in $\epsilon$, formally constant with respect
to the dynamics. Under relatively mild assumptions, locally there exist actual
functions which are invariant to all orders in $\epsilon$ and, under further
conditions, even within exponentially small errors (see for instance
\cite{RamisS} and also the literature cited there).

In the present paper we show that under suitable analyticity conditions,
adiabatic invariants are Borel summable and their Borel sums are {\em exact
  invariants} in the regions of regularity.

\centerline{*}

The technique that we use applies to a wide class of differential
systems but for the sake of clarity we prefer to focus on
relatively simple equations and discuss later how the results and
methods extend. Our prototypical equation is

\begin{equation}\label{1}
\psi'' +(\epsilon^{-2}-V(x,\epsilon))\psi = 0\qquad (\epsilon\to
0),
\end{equation}

\z a special case of which is the {\em one dimensional Schr\"odinger
  equation} in the large energy limit with analyticity and decay
conditions in some strip in $\CC$.

A number of different equations in which formal invariants arise can be
easily
brought precisely to the form (\ref{1}).

The {\em parametrically perturbed
  pendulum}
\begin{equation}\label{pend}
\ddot x + \omega(t\epsilon)^2\;x=0.\qquad (\epsilon\to 0)
\end{equation}
\z after the substitution
$x(t)=\omega^{-1/2}(\epsilon t)f(\int_a^{\epsilon t}\omega(s)ds)$
becomes

$$f''+\left(\epsilon^{-2}+\frac{3}{4}\frac{{\omega'}^2}{\omega^4}-\frac{\omega''}{2\omega^3}\right)f=0.$$
Our results do not apply as such to equations with periodic
coefficients but, by transformations, some of these equations can be
brought to our setting.

\z For example, {\em Mathieu's equation (\cite{Kamke} p. 404)-- in a singular perturbation regime} 
$$
\epsilon^2\ddot h - (a\cos(2x)+b)h=0
$$

\z can be transformed by taking
$h=\sin^{-1/2}(2x)\,f(\alpha\cos(2x))$, $\sigma=\alpha\cos(2x)$, to

$$\epsilon^2
f''+\frac{1}{4}\left(\frac{2a\sigma+b\alpha}{\alpha(\alpha^2-\sigma^2)}+
  \epsilon^2\frac{2\alpha^2+\sigma^2}{(\alpha^2-\sigma^2)^2}\right)f=0$$

\z 

Also, as it will become transparent, the equations for which the
methods in this paper apply can be higher order and/or contain
nonlinear terms.

\centerline{*}

\z We let
\begin{equation}
  \label{eq:deflam}
  \lambda=\epsilon^{-1}.
\end{equation}
Eq. (\ref{1}) admits the formal solutions\footnote {This type of
  expansion is best suited for Borel summability in our setting.}

$$A_+(\lambda)e^{i\lambda
x}\sum_{k=0}^{\infty}g_{k,+}(x)\lambda^{-k}+A_-(\lambda)e^{-i\lambda x}\sum_{k=0}^{\infty}g_{k,-}(x)\lambda^{-k}$$

\z with

$$g_{0,\pm}=1, g_{1\pm}=\mp\frac{1}{2} i\int V(s)ds,...$$

We consider (\ref{1}) on a bounded open interval
$I\subset\mathbb{R}$ with initial condition prescribed at some
$x_0\in I$. Without loss of generality we may take $I=(-1,1)$ and
$x_0=0$.

We also assume that
\[\psi(0)=f(\lambda),\]
where $f$ is the sum of a convergent or Borel-summable series.

Throughout this paper, Borel summation is understood in the following
way:
\begin{Definition}
A Borel-summable series $\tilde y:=\sum_{k=K}^\infty y_k\la^{-k}$,
$K\in\ZZ$ is a formal
power series with the following properties
\begin{enumerate}
\item the truncated Borel transform $Y=\mathcal{B}\tilde
y:=\sum_{k>0}\frac{y_{k}}{(k-1)!}t^{k-1}$ of $\tilde y$ has a
nonzero radius of convergence, \item $Y$ can be analytically
continued along $[0,+\infty)$ and \item the analytic continuation
$Y$ grows at most exponentially along $[0,+\infty)$ and is
therefore Laplace transformable along $[0,+\infty)$.
\end{enumerate}
The Borel sum $y$ of $\tilde y$ is then given by
\begin{equation}
  \label{eq:defLB}
  y = \lap\bor\tilde{y}:=\sum_{k=K}^{0}y_{k}\la^{-k} +
  \mathcal{L}Y,
\end{equation}
where the sum is understood to be zero if $K>0$ and $\mathcal{L}$
denotes the usual Laplace transform.
\end{Definition}
\begin{Remark}
  We note that although we require conditions in a $\CC-$ neighborhood
  $\mathcal{N}$ of $[0,\infty)$ rather than in a sector, uniqueness of
  $y$ follows since $Y$ is analytic in $\mathcal{N}$ and is uniquely
  defined near zero by $\mathcal{B}\tilde y$.
\end{Remark}

\z Our main results are the following :
\begin{Theorem}\label{T1}
Assume that $V$ is analytic in $\epsilon$ and $x$:

\begin{equation}\label{g1}V(x,\epsilon)=\sum_{k=0}^{\infty}V_k(x)\epsilon^k\end{equation}

\z with $V_k$ analytic in the strip $\mathcal{S}=\{-1<\Re(x)<1\}$,
real valued on the real line and, for some $B,K,\delta>0$,
satisfying
\begin{equation}\label{growth}
\sup_{k\in\NN,x\in\mathcal{S}}(1+|x|^{1+\delta})B^{-k}|V_k(x)|=K<\infty.
\end{equation}
Then given $\la=\epsilon^{-1}>0$
%(cf. \ref{eq:deflam})
large enough,
the general solution of (\ref{1}) in $(-1,1)$ can be written in
the form
$$
\psi = C_1 e^{i\la x}\phi_{+}(x;\la) + C_2 e^{-i\la x}\phi_{-}(x;\la)
$$
where $\phi_+$ and $\phi_-$ are conjugate expressions of each other
and are the Borel sums of their asymptotic power series ($\sum
a_{k_+}(x)\la^{-k}$ and $\sum a_{k_-}(x)\la^{-k}$ respectively). \ \z
The initial condition is taken to be
\begin{equation}\label{IC}
\psi(0)=f(\lambda)=1/\la.
\end{equation}
Furthermore, $\phi_+$ (and $\phi_-$) are uniquely determined in
the following sense : if for some neighbourhood $\mathcal{V}$ of
$0$ and some $\nu>0$,
$$
\phi(x;\la)=\mathcal{L}(\chi(x,\cdot)),
$$
where $\chi(\cdot,t)\in C^2(\mathcal{V})$ for $t\in (0,\infty)$
and $\chi(x,\cdot)\in L^1_\nu:=L^1(\mathbb{ R}^+,e^{-\nu t}dt)$ for
$x\in(-1,1)$,
and if $e^{i\la x}\phi$ solves (\ref{1}),(\ref{IC}) on some
neighbourhood of $0$ then
$$
\phi=\phi_+.
$$
\end{Theorem}

\begin{Theorem}\label{T2} Consider a $V$ as in
  Theorem \ref{T1} and take $\psi$ to be a solution of (\ref{1}),(\ref{IC}) in
  a neighbourhood of $0$, which can therefore be written as
$$
\psi = C_1 e^{i\la x}\phi_{+}(x;\la) + C_2 e^{-i\la
x}\phi_{-}(x;\la).
$$
Then, in the region where the assumptions are satisfied,
$C(x;\lambda;\psi,\psi'):=C_1C_2$ is an {\bf exact invariant} of
(\ref{1}) and is the Borel sum of an {\bf adiabatic invariant}
\begin{eqnarray}
  \label{eq:definvar}\tilde{C}=\sum
c_k(x;\psi,\psi')\la^{-k},\qquad\text{i.e.}
\end{eqnarray}
$$\ C(x;\lambda;\psi,\psi')=\lap\bor \tilde{C}.$$

\end{Theorem}
For example, in (\ref{pend}), in the regions where the analyticity
assumptions are fulfilled there exists an {\em actual} invariant
of the form $C(t;\epsilon;x,\dot x)$. $C$ is analytic in a sector
in $\epsilon>0$, {\em
  analyzable} in the given sector at $\epsilon=0$ (in this case, it
simply means that the Taylor series of the function is Borel
summable to the function itself), and it is straightforward to
check that to leading order it assumes the familiar form

$$C\sim E/\left(2\omega^2\right);\ \ E=\frac12 (\dot{x}^2+\omega^2 x^2)$$

\begin{Remark}
All of our results still hold if instead of \eqref{IC}, one
assumes only that $\psi(0) = f(\lambda)$, where $f$ is the sum of
a Borel summable series. This is a simple consequence of the fact
that equation \eqref{1} is linear and that Borel summable series
form an algebra (see Lemma \ref{L2}).

\end{Remark}

\begin{Remark}All of our results also hold if instead of \eqref{growth}, the coefficients $V_k$ satisfy the following weaker
growth assumption at infinity:
\[|V_k(x)|\le B^k g(|x|)\qquad\text{for some $B>0$ and $g\in L^1(\RR)$}.\]
\end{Remark}

\section{Proof of Theorem \ref{T1}}
\z The core of the proof is contained in the case where $V$ is independent of
$\epsilon$ and we start with this case. In the following $\mathcal{C}$ denotes
any constant the value of which is not significant to the analysis.

\subsection{Formal derivation of a fixed-point equation}
First, we (formally) manipulate (\ref{1}) to obtain a fixed point
problem.
To leading order, the solutions of (\ref{1}) are $e^{i\la x},
e^{-i\la x}$ so that we look for solutions of the form
\begin{equation}\label{eqphi+}
\psi=\psi_+= e^{i\la x}\phi_+(x)\qquad\text{and}\qquad
\psi=\psi_-= e^{-i\la x}\phi_-(x).
\end{equation}
Then $\phi_+$ solves
\begin{equation}\label{2}
\phi'' + 2i\la \phi' - V(x)\phi=0.
\end{equation}

We then seek for solutions of (\ref{2}) of the form
$\phi(x;\la):=\mathcal{L}(\chi(x,\cdot))$, where $\mathcal L$ is
the usual Laplace transform. One can easily check that
$$
\frac{\mathcal{L}(\chi(x,\cdot))}{\lambda}=\mathcal{L}(\Psi(x,\cdot))\qquad\text{where}\qquad
\Psi(x,t)=\int_0^t \chi(x,\tau)\;d\tau,
$$
so that dividing (\ref{2}) by $\lambda$ leads to
\begin{equation}\label{eqPsi}
\left\{
\aligned
\Psi_{xx} + 2i \Psi_{xt}  = &V(x)\Psi \\
\Psi(x,0)=&0\\
\Psi(0,t)=&\int_0^t d\tau = t
\endaligned
\right.
\end{equation}
To control this equation, we pass to bicharacteristic coordinates,
namely we let
\begin{equation}
s=-2ix+t
\label{bichar}
\end{equation}
and obtain for $\Phi(s,t):=\Psi(x,t)$,
\begin{equation}\label{eqphi}
\left\{
\aligned
\Phi_{st}\ \ \ \  & = \frac{1}4 V(i(s-t)/2)\Phi \\
\Phi(s,0)&=0\\
\Phi(t,t)&=t
\endaligned
\right.
\end{equation}

\z Integrating (\ref{eqphi}) yields
\begin{equation}\label{integraleq}
\Phi(s,t)= J(\Phi)(s,t),
\end{equation}

\z where

\begin{equation}
  \begin{aligned}
    J(\Phi)(s,t)=t-\frac{1}{4}\int_0^t\int_{t}^sV(i(s_1-t_1)/2)
\Phi(s_1,t_1)d s_1d t_1:=&\\
t+t(s-t)\int_0^1\int_0^1V\big[i((1-\alpha)s+(\alpha-\beta)t)/2\big]\Phi\big[(1-\alpha)s+\alpha
t,\beta t\big]\;d\alpha\; d\beta.
  \end{aligned}
  \label{eq:defL}
\end{equation}

\subsection{Solving the fixed-point equation: analyticity and
exponential bounds}
It is useful to note the relation between $\Psi$ and $\Phi$:
  $$\Psi(x,t)=\Phi(-2ix+t,t)\ \ \mathrm{or}\ \ \
  \Phi(s,t)=\Psi(i(s-t)/2,t).$$

For $\nu>0$ and $\Omega$ an open set of $\mathbb{C}^2$, let now
$\mathfrak{B}_\nu(\Omega)$ be the set of analytic functions $F$ over $\Omega$
with finite norm
\begin{equation}
  \label{eq:norm1}\|F\|_\nu=\sup_\Omega |F(s,t)|e^{-\nu|t|}.
\end{equation}
Let
$$P=\{(s,t):|s|<R,|t|<R\} $$

\z and

\[S_1=\{(s,t): \Re(s),\Re(t)>0,\;\
\max(|\Im(t)|,|\Im(s)|,|\Im(t-s)|)<2\}\cup P.\]

\begin{Lemma}\label{L1}\ Assume $V(x)$ is analytic in $S=\{|x|<R\}\cup
\{x:
%\Im(x)>0;\
|\Re(x)|< 1\}$ and
  $|V(x)|<K (1+|x|)^{-1-\delta}$ in $S$ for some $\delta>0$.
\z Then
\[J:\mathfrak{B}_\nu(S_1)\to\mathfrak{B}_\nu(S_1)\]
is a contraction for $\nu>0$ large enough.

\end{Lemma}

\begin{proof}
First we check that $J$ is well-defined on $\mathfrak{B}_\nu(S_1)$.
Clearly,
 $$(s,t)\in P\Longrightarrow i((1-\alpha)s+(\alpha-\beta)t)/2\in
\{x:|x|<R\}$$
 for $\alpha,\beta\in[0,1]$.
%Now, if $(s,t)\in S_1\setminus P$ then
%\[
%\Im(i((1-\alpha)s+(\alpha-\beta)t)/2)=\Re((1-\alpha)s+(\alpha-\beta)t)/2>\Re((1-\beta)t)/2>0.
%\]
Also
\[
|\Re(i((1-\alpha)s+(\alpha-\beta)t)/2)| = |(1-\alpha)\Im(s) +
(\alpha-\beta)\Im(t)|/2<1,
\]
provided this inequality holds for
$(\alpha,\beta)=(0,0),(0,1),(1,0),(1,1)$ since the above expression is linear in $(\alpha,\beta)$. The
condition is thus

$$\max(|\Im(t)|,|\Im(s)|,|\Im(t-s)|)<2,\qquad\text{i.e. $(s,t)\in
S_1$.}$$
It follows easily that $J$ maps analytic-over-$S_1$ functions to
themselves.
Also, $\|t\|_\nu\le \nu^{-1}$. We are left with proving that $J$ is
contractive
(in a nonlinear case, one would restrict $J$ to a ball of radius
$\mathcal{C}\nu^{-1}$) and this follows from the following estimate
\begin{multline}
  \label{eq:evalint}
\left|  \int_s^t V(i(s_1-t_1)/2)
F(s_1,t_1)ds_1\right|\le\\
\mathcal C\|F\|_\nu e^{\nu
|t_1|}\int_{-\infty}^{\infty}(1+|u|)^{-1-\delta}du
= \mathcal C\|F\|_\nu e^{\nu |t_1|},
\end{multline}
\z whence by integration
$$\|J(F)\|_{\nu}\le \mathcal{C}\nu^{-1}\| F\|_{\nu}.$$
\end{proof}
Applying Lemma \ref{L1} with $R=1$ and Picard's fixed-point
theorem, we thus obtain a solution $\Phi\in\mathfrak{B}_\nu(S_1)$
of \eqref{integraleq}.
\subsection{Solving the original equation}
\begin{proof}
This part of the proof is standard but we include it for convenience.
To go back to the original equation, we just need to reverse and
justify the transformations that lead us from (\ref{1}) to
(\ref{integraleq}) : since $\Phi\in\mathfrak{B}_\nu(S_1)$, $\Psi$ is
analytic in a
neighbourhood of $(-1,1)\times(0,\infty)$ and
$$
|\Psi(x,t)|\leq \mathcal Ce^{\nu |t|}.
$$
Clearly, $\Psi$ satisfies (\ref{eqPsi}). We claim that given any
$x_0\in(-1,1)$ there exists $r>0$ such that each partial
derivative of $\Psi$ is exponentially bounded on
$(x_0-r,x_0+r)\times(0,\infty)$ so that $\phi=\la\mathcal
L(\Psi(x,\cdot))$ is well-defined and will solve (\ref{2}) with
the initial condition $\phi(0;\la)=1/\la$. Indeed,
$$
\begin{matrix}
\Psi_x=-2i\Phi_s,&\Psi_{xx}=-4\Phi_{ss}\\
\Psi_t=\Phi_s+\Phi_t,& \Psi_{xt}=-2i(\Phi_{ss}+\Phi_{st})
\end{matrix}
$$
But since $\Phi$ is analytic in $U=P\cup S_1$, letting
$s_0=-2ix+t_0$ for $x\in(x_0-r,x_0+r)$ and $t_0\in\RR^+$, we have
that $(s_0,t_0)\in U$ and by Cauchy's formula,
$$
|\Phi_s(s_0,t_0)|=(2\pi)^{-1}\left|\int_{\mathcal
S(s_0,r)}\Phi(s,t_0)/(s-s_0)^2\; ds\right| \leq \mathcal C_r
e^{\nu |t_0|},
$$
where $\mathcal S(s_0,r)$ denotes the circle centered at $s_0$ of
radius $r=(1-|x_0|)/2$. The bounds on the other derivatives can be
obtained similarily and we therefore take Laplace transforms in
(\ref{eqPsi}) to obtain a solution of (\ref{1}),(\ref{IC}) on
$(-1,1)$, of the form
$$
\psi(x;\la)=\psi_+=e^{i\la x}\phi_+=\la e^{i\la x}\mathcal
L(\Psi(x,\cdot))\qquad\text{where $\qquad|\Psi(x,t)|\leq
\mathcal Ce^{\nu|t|}$.}
$$
Working with $\phi_-=\overline{\phi_+}$ (as defined in (\ref{eqphi+})),
we
obtain similarily a solution of the form $\psi(x;\la)=e^{-i\la
x}\phi_-$.

In the following we will rely on the fact that Borel summable
series are closed under algebraic operations.

\begin{Lemma}\label{L2}
Borel summable series form a field.
\end{Lemma}
\z Though rather straightforward (see \cite{Balser}, \cite{Ramis},
\cite{R-M}, \cite{R-S}), we provide a proof for convenience of the
reader.  The fact that Borel summable series form an {\em algebra}
is shown in \cite{Balser}.  Given a series $\tilde y$, we want to
construct its multiplicative inverse $(\tilde y)^{-1}$. Up to
factoring out a monomial $y_K\la^{K-1}$ in the expension of
$\tilde y$, we may always assume that
$$
\tilde y= \lambda(1 + \tilde f)\qquad\text{for some Borel summable
$\tilde f=o(1)$}.
$$
The inverse $g$ of $\tilde y$ must then satisfy
\begin{equation}
(1+\tilde f)\tilde g = 1/\la.
\label{eq:oneplusfg}
\end{equation}
If $F=\mathcal{B}\tilde f$ is the Borel transform of $f$, defined
on a fixed (star-shaped-about-the-origin) neighbourhood $\Omega$
of $[0,\infty)$ and if $G$ denotes the Borel transform of the
inverse $g$ we are looking for, we must have
\begin{equation}\label{convolution} G = 1 - F*G, \end{equation}
where
\begin{equation}\label{defconvolution}
(F*G)(t)=\int_0^t F(q)G(t-q)dq:= t\int_0^1 F(\alpha
t)G((1-\alpha)t)d\alpha.
\end{equation}
Given $A>0$ and $\nu>0$,consider the norm
\begin{equation}\label{norma}
\|F\| := \sup_{t\in\Omega} (A+|t|)^2|F(t)|e^{-\nu t},
\end{equation}
defined for $F\in\mathcal{B}_{\nu,A}(\Omega)$, the space of
(exponentially bounded) analytic functions over $\Omega$, equipped
with the above norm. \eqref{convolution} will have a solution $G$
provided $G\to F*G$ is a contraction in that space. Now,
\begin{multline}\label{convestimate}
|F*G|\le \|F\| \|G\| e^{\nu|t|}\int_0^{|t|}
(A+u)^{-2}(A+|t|-u)^{-2}du \\= \|F\| \|G\|
e^{\nu|t|}\left(\int_0^{|t|/2} + \int_{|t|/2}^{|t|}\right)
(A+u)^{-2}(A+|t|-u)^{-2}\;du
\end{multline}
and
\begin{equation*}
\int_0^{|t|/2} (A+u)^{-2}(A+|t|-u)^{-2}du\le
(A+|t|/2)^{-2}\int_0^\infty (A+u)^{-2}du \le
\frac{\mathcal{C}}{A}(A+|t|)^{-2}.
\end{equation*}
Working similarly with the last integral in \eqref{convestimate},
we obtain that
\begin{equation}\label{convnorm}
\|F*G\|\le \frac{\mathcal{C}}{A}\|F\| \|G\|.
\end{equation}
Hence, fixing $A>0$ large enough it follows that $G\to F*G$ is
contractive.

Taking Laplace transforms and applying Watson's lemma \cite{Orszag}, we obtain
the (unique) solution $\tilde g$ of (\ref{eq:oneplusfg}) defined
by
$$
\mathcal{L}G\sim\tilde g.
$$
\QED

\subsection{Uniqueness of $\phi_+$}

Take $\phi=\mathcal L(\chi(x,\cdot))$ and $\mathcal{V}$ a neighbourhood of
the origin as in the statement of Theorem \ref{T1}. As in Step 2.1,
using the same notations to go from $\phi$ to
$\Phi$, it follows that $\Phi$ solves (\ref{integraleq}), so
that we only need to prove that $J$ is a contraction for the norm
$$
\|F\|_Y = \sup_x\int_0^\infty |F(x,t)|e^{-\nu t}dt
$$
and observe that since $L^1_\nu\subset L^1_{\nu'}$ when $\nu<\nu'$, we
may choose
$\nu>0$ as large as we please.
Rewriting (\ref{eq:defL}) in terms of $\chi\in Y$, we have
\begin{multline}
K(\Phi)(s,t):= -4(J(\Phi)(s,t) -t)\\
  = t(s-t)\int_0^1\int_0^1\int_0^{\beta t} V(x_{\alpha\beta})
\chi(x_{\alpha\beta},\tau)d\tau\; d\alpha\; d\beta,
\end{multline}
where
$$
x_{\alpha\beta} = \frac{i}{2}\left((1-\alpha)s+(\alpha-\beta)t\right).
$$
So
\begin{equation}
|K(\Phi)|\le \|\chi\|_Y |t(s-t)|\int_0^1\int_0^1
|V(x_{\alpha\beta})| e^{\nu\beta t}d\alpha d\beta \le
\mathcal{C}\nu^{-1}\|\chi\|_{Y},
\end{equation}
where we used (\ref{eq:evalint}) in the last inequality.
\label{eq:uniqueness}
%\left|\int_0^t\int_{s}^tV(i(s_1-\tau)/2)F(s_1,\tau)d s_1d\tau\right|
%\\\le (|s|+|t|)v\|F\| \int_0^{|t|} e^{\nu | \tau|}d|\tau|\le
%\frac{2Rv\|F\|}{\nu} e^{\nu |% t|}=:\mathcal{C}{\nu}^{-1}\|F\|e^{\nu|t|}

\subsection{$\psi_+$ and $\psi_-$ are linearily independent}
Observe that $\psi_+$ and $\psi_-$ defined in (\ref{eqphi+}) are
conjugate expressions of each other and
arguing by contradiction, suppose that for some $A=A(\la)\in S^1$
$$
\psi_+=A \psi_-
$$
Let $B(x;\la):=\phi_+/\phi_-$. Then
\begin{equation}
B = Ae^{-2i\la x}.
\label{eqB}
\end{equation}

By the Lemma \ref{L2}, $B$ defined by (\ref{eqB}) is asymptotic to a Borel summable power
series and since
$\phi_+(0;\la)=\phi_-(0;\la)=1/\la$, we have for $x=0$,
\begin{equation}
A(\la)=B(0;\la)=1+o(1) \qquad\text{as $\la\to\infty$,}
\label{A=B}
\end{equation}
whereas for any other fixed value of $x\neq0$, there exist
$K=K(x)\in\ZZ$ and $b=b(x)\in\mathbb{C}$ such that
$$
A(\la)e^{-2i\la x}=B= \la^K(b+o(1)).
$$
Combining this equation with (\ref{A=B}), we obtain
$$
(1+o(1))e^{-2i\la x} =  \la^K(b+o(1)),
$$
which is clearly impossible.

\subsection{The case (\ref{g1})} We explain here how to adapt
the argument for a potential $V=V(x,\e)$ depending on
$\e=1/\lambda$. According to \eqref{g1}, if $\phi(x)=\mathcal
L(\chi(x,\cdot))$, the inverse Laplace transform of $V\phi$ is
given by:
\begin{equation}
  \label{eq:il2}
\mathcal{L}^{-1}V\phi=
V_0(x)\chi+\sum_{k=1}^{\infty}\frac{V_k(x)}{(k-1)!}\;t^{k-1}*\chi,
\end{equation}
where $*$ is the convolution product defined in
\eqref{defconvolution}. Hence, taking as before
$\Psi(x,t)=\int_0^t\chi(x,\tau)\;d\tau$, we obtain
$$
\frac{V\phi}{\lambda}= \mathcal{L}\left(V_0(x)\Psi +
\sum_{k=1}^\infty\frac{V_k(x)}{(k-1)!}\;t^{k-1}*\Psi\right).
$$
Following the analysis of the $\epsilon$-independent case, we end
up with the following operator :
\begin{equation*}
    \tilde J(\Phi)(s,t)=J(\Phi)(s,t)-\frac{1}{4}\sum_{k>0}\int_0^t\int_{t}^s
\frac{V_k(i(s_1-t_1)/2)}{(k-1)!}[t^{k-1}*\Phi](s_1,t_1)\;d s_1\;d
t_1,
\end{equation*}
where $J$ is defined by\eqref{eq:defL}. Again, we must prove that
$\tilde J$ is contractive. To do so, instead of (\ref{eq:norm1}),
we use the norm $\|\cdot\|$ defined in \eqref{norma}. For $A$
large enough, it follows from \eqref{convnorm} that
\begin{equation}\label{convnorm2}
\|F*G\|\le\|F\|\|G\|\text{ holds for all $F$, $G$ analytic in
$\Omega$.}
\end{equation}
In order to control $\|\tilde J\|$, we perform a few side
computations. We first estimate $\|t^{k-1}\|$ :
$$
\|t^{k-1}\|\le \mathcal{C}\sup_{t\in\RR^+}\left({t^{k+1}e^{-\nu
t}}+{t^{k-1}e^{-\nu
t}}\right)=\mathcal{C}\left[\left(\frac{k+1}{\nu}\right)^{k+1}e^{-(k+1)}+\left(\frac{k-1}{\nu}\right)^{k-1}e^{-(k-1)}\right].
$$
By Stirling's formula, it follows that
\begin{equation}\label{stirling}
\|t^{k-1}\|\le \mathcal C (k-1)!\;\;k^{3/2}\nu^{-k+1}.
\end{equation}
Next, we estimate for $t>0$ the following quantity :
$$
e^{-\nu t}(A+t)^2 \int_0^t e^{\nu t_1}
(A+t_1)^{-2}dt_1=\left(\int_0^{t/2}+\int_{t/2}^t\right)e^{\nu(t_1-t)}\left(\frac{A+t}{A+t_1}\right)^2\;dt_1.
$$
On the one hand
$$
\int_0^{t/2}e^{\nu(t_1-t)}\left(\frac{A+t}{A+t_1}\right)^2\;dt_1\le
e^{-\nu t/2}\int_0^{t/2}\left(\frac{A+t}{A+t_1}\right)^2\;dt_1\le
e^{-\nu t/2}\frac{(A+t)^2}{A}\le\frac{\mathcal C}{\nu^2},
$$
while on the other hand
$$
\int_{t/2}^{t}e^{\nu(t_1-t)}\left(\frac{A+t}{A+t_1}\right)^2\;dt_1\le
\mathcal C \int_{t/2}^t e^{\nu(t_1-t)}\;dt_1\le \mathcal C/\nu.
$$
Hence,
\begin{equation}\label{estimint}
e^{-\nu t}(A+t)^2 \int_0^t e^{\nu t_1} (A+t_1)^{-2}dt_1\le\mathcal
C/\nu.
\end{equation}
We can now estimate $\|\tilde J(\Phi)\|$. Instead of
\eqref{eq:evalint}, we have
\begin{multline*}
\left|  \int_s^t V_k(i(s_1-t_1)/2)
[t^{k-1}*\Phi](s_1,t_1)ds_1\right|\le \mathcal C\|t^{k-1}*\Phi\|
e^{\nu
|t_1|}(A+|t_1|)^{-2}B^k\int_{-\infty}^{\infty}(1+|u|)^{-1-\delta}du
\\ \le \mathcal C\nu \left(\frac B\nu\right)^k (k-1)!\;\;k^{3/2} e^{\nu
|t_1|}(A+|t_1|)^{-2}\|\Phi\|,
\end{multline*}
where we used \eqref{convnorm2} and \eqref{stirling} in the last
inequality. Whence by integration and by \eqref{estimint}, we
finally obtain
$$
\|\tilde J(\Phi)\|\le
\mathcal{C}\left(\nu^{-1}+\sum_{k>0}k^{3/2}\left(\frac
B\nu\right)^k\right)\|\Phi\|\le\frac{\mathcal C}{\nu}\|\Phi\|.
$$
\end{proof}
\subsection{Further generalisations} It is not difficult to modify the
proof given to allow for higher order equations, which, possibly after
transformations, have relatively simple bicharacteristics, or for
nonlinear dependence on $\psi$ which does not affect the highest
derivative. In the nonlinear case, the initial condition must
obviously be left in the form of a general sum of a convergent or
Borel summable series.  The strategies of dealing with nonlinearities
are described in \cite{Duke} and \cite{CT}.

\section{Proof of Theorem \ref{T2}} This follows easily from the
explicit expression
of $C$, Theorem~\ref{T1} and Lemma~\ref{L2}.
Indeed, differentiating $\psi$ (with respect to $x$), we obtain
$$
\aligned
\psi&=C_1e^{i\la x}\phi_+ + C_2e^{-i\la x}\phi_-\\
\psi'&=C_1e^{i\la x}A + C_2e^{-i\la x}B
\endaligned
$$
where
$$
A=\phi_+' +i\la\phi_+\qquad\text{and}\qquad B=\phi_-' -i\la\phi_-
$$
So that
$$
C:=C_1C_2=
\dfrac{
\left|
\begin{matrix}\psi  & \phi_-\\
        \psi' & B
\end{matrix}
\right|
\cdot\left|
\begin{matrix}\phi_+ & \psi\\
             A & \psi'
\end{matrix}
\right|
}
{
{\left|
\begin{matrix}\phi_+ & \phi_-\\
             A & B
\end{matrix}
\right|}^2
}
$$
By construction, $\phi_+\sim\sum a_k(x)\la^{-k}$ and is uniquely
determined so that using Lemma~\ref{L2}, we may conclude that
$$
C\sim\sum c_k(x,\psi,\psi')\la^{-k}
$$
where all considered power series are Borel summable.
$\Box$

\section{Discussion of singularities of $\Psi$}
So far we have assumed that $V$ was analytic. If this is not so,
singularities of $\Psi$ can originate in the singular points of
$V$.  We restrict the analysis to the relatively common situation
where $V$ has a branch-point at some point $x_0$ (we may  assume
without loss of generality $x_0=0$) of the form
$$V(x)=x^{-\beta}V_1(x),$$
\z with $V_1$ analytic at 0. We assume
$0<\Re \beta<1$ (for different $\Re \beta$, the analysis can be done
similarly).

 We let $P=\{(s,t):|s|<\frac{1}{2\tau},|t|<\tau\}$ and
consider the region

\[S_1=\{(s,t): \Re(s)>\Re(t)>0,\;\
\max(|\Im(t)|,|\Im(s)|,|\Im(t-s)|)<2\}\cap P\]

\z Assume that $(s,t)\to V_1(i(s-t)/2))$ is analytic in $S_1$ and
that $V$ satisfies the following estimate (similar to that in
Lemma \ref{L1}) : for some $\delta>0$, $K>0$,
$$|V(x)|\le K (1+|x|)^{-1-\delta}\quad\text{ if }\quad |x|>1.$$
\begin{Proposition}\label{Pr.6}
   There exists a solution $\Phi\in C^2(S_1)$ of \eqref{eq:defL}. Furthermore, letting $v_0=V_1(0)$, $\Phi$ satisfies

\begin{multline}\label{struct}\Phi(s,t)=t+\frac{v_0}{4}\Big\{ s\frac{
    s^{2-\beta}-(s-t)^{2-\beta}}{(1-\beta)(2-\beta)(3-\beta)}+\frac{(s-t)^{2-\beta}}{(1-\beta)(3-\beta)}t\\+\frac{
    t^{3-\beta}}{(1-\beta)(2-\beta)(3-\beta)}\Big\}+O(t^3),\
   \end{multline}
as $t\to 0.$
\end{Proposition}
\z {\bf Remarks}

(1) In particular, it follows easily that for $s\ne 0$, $\Phi$
does not extend analytically at $t=0$. Indeed, assuming the
contrary, we would have $Ct^{3-\beta}=A(t)+O(t^3)$ with $A$
analytic; as $t\to 0$ this forces $A(t)\sim Ct^{3-\beta}$. But $A$
is analytic and we must have for some $n\in\NN$ and $C'\in\RR$,
$A(t)\sim C't^n$, which is a contradiction.

(2) Proposition \ref{Pr.6} provides information about the analytic
continuation (still denoted $\Phi$) of the solution to
(\ref{eq:defL}). Indeed, we can work as in Lemma \ref{L1} to show
analyticity in a region where $t>\tau>0$. Lemma \ref{L1} and
Proposition \ref{Pr.6} provide uniqueness in the space of
exponentially bounded analytic functions over the corresponding
region, so that by obvious imbeddings the analytic continuation
$\Phi$ coincides with the fixed point of \eqref{eq:defL}. Changing
variables, i.e. going back to $\Psi$, we conclude that $\Psi$ is
nonanalytic at $t=0$. This implies in turn a Stokes transition on
$\phi$ (and thus the adiabatic constant constructed in Theorem 3).

(3) Furthermore, the fixed point procedure that leads to relation
(\ref{struct}) can easily (and rigorously) provide more detailed
information about the singularity manifold; it suffices to
construct the space of functions in which the fixed point equation
(\ref{eq:eqH}) below is considered in such a way that the norms
entail information about the singularity type to be proved: the
details are quite straightforward but at the same time quite long
and we will not elaborate on them in the present paper. Similar
constructions can be found in \cite{Duke}.

\bigskip

 \z {\em Proof.} From (\ref{eq:defL}) we have

\begin{equation}
  \label{eq:eqH}
  \Phi(s,t)=t-\frac{1}{4}\int_0^t\int_{t}^s
\frac{V_1(i(s_1-t_1))}{(s_1-t_1)^\beta}
\Phi(s_1,t_1)d s_1dt_1
\end{equation}

\z The fact that $J$ is defined and contractive with norm $O(t)$ on the
functions defined in $S_1$ with the sup norm, follow as in Lemma \ref{L1}.

We have $(1-J)\Phi=t$ and thus for small $t$, $\Phi=t+Jt+O(J^2t)=t+Jt+O(t^3)$
and, with $v_0=V_1(0)$ the result follows.\QED

\bigskip

\z {\bf Acknowledgments}.  Work partially supported by NSF Grants
0103807 and 0100495. O C would like to thank Prof. G Hagedorn and A
Joye for interesting discussions and comments.

\end{document}